\documentclass{amsart}
\usepackage{amsmath}
\usepackage{amsfonts}
\usepackage{amssymb}
\usepackage[curve,matrix,arrow]{xy}
\def\qedsymbol{~\hfill QED.}
\def\proofname{Proof.}
\newenvironment{Proof}{\par\noindent{\it\proofname}}{{\unskip\nobreak\hfill{\it\qedsymbol}}\par\vskip 9pt}
\newenvironment{Proof*}{\par\noindent}{{\unskip\nobreak\hfill{\it\qedsymbol}}\par\vskip 9pt}
\ifx\undefined\bysame \newcommand{\bysame}{\leavevmode\hbox to3em{\hrulefill}\,}\fi
\newtheorem{Thm}{Theorem}[section]
\newtheorem{Lem}[Thm]{Lemma}
\newtheorem{Prop}[Thm]{Proposition}
\newtheorem{Fact}[Thm]{Fact}
\newtheorem{Cor}{Corollary}[Thm]
\newtheorem{Conj}[Thm]{Conjecture}

\newtheorem{Rem}[Thm]{Remark}

\newcommand{\wcat}{\operatorname{{\it w}cat}}

\newcommand{\cat}{\operatorname{cat}}
\newcommand{\Cat}{\operatorname{Cat}}

\newcommand{\integral}{\mathbb Z}
\newcommand{\complex}{\mathbb C}

\newcommand{\homeo}{\approx}

\newcommand{\comp}{\smash{\lower-.1ex\hbox{\scriptsize$\circ$}}}
\newcommand{\fatvee}{\operatorname{T}}
\newcommand{\widebar}[1]{\overline{#1}}
\title[L-S categories of Lie groups]{L-S categories of simply-connected compact simple Lie groups of low rank}
\author[N.~Iwase]{Norio Iwase}
\address[N.~Iwase]{Faculty of Mathematics, Kyushu University, Ropponmatsu Fukuoka 810-8560, JAPAN.}
\email{iwase@math.kyushu-u.ac.jp}
\author[M.~Mimura]{Mamoru Mimura}
\address[M.~Mimura]{Department of Mathematics, Okayama University, Okayama 700-8530, JAPAN.}
\email{mimura@math.okayama-u.ac.jp}
\thanks{The second named author is partially supported
by the Grant-in-Aid for Scientific Research \#12640025 from the Japan Society of Promotion of Science.}
\begin{document}
%
%
\begin{abstract}
We determine the L-S category of $Sp(3)$ by showing that the $5$-fold reduced diagonal $\widebar{\Delta}_5$ is given by $\nu^2$, using a Toda bracket and a generalised cohomology theory $h^*$ given by $h^*(X,A) = \{X/A,{\mathbb S}[0,2]\}$, where ${\mathbb S}[0,2]$ is the $3$-stage Postnikov piece of the sphere spectrum ${\mathbb S}$.
This method also yields a general result that $\cat(Sp(n)) \geq n+2$ for $n \geq 3$, which improves the result of Singhof \cite{Singhof:ls-lie-group-ii}.
\end{abstract}
%
%
\date{August 29th, 1995}
\date{October 22nd, 2001, First draft}
\date{November 26th, 2001, Second draft}
\date{December 1st, 2001, third draft}
\date{December 17th, 2001, third draft}
%
%
\maketitle
\section{Introduction}
\par\noindent
{\sc In this} paper, each space is assumed to have the homotopy type of a CW complex.
The (normalised) L-S category of $X$ is the least number $m$ such that there is a covering of $X$ by $(m+1)$ open subsets each of which is contractible in $X$.
Hence $\cat{\{\ast\}} = 0$.
By Lusternik and Schnirelmann \cite{LS:l-s_category}, the number of critical points of a smooth function on a manifold $M$ is bounded below by $\cat{M}+1$.

G. Whitehead showed that $\cat(X)$ coincides with the least number $m$ such that the diagonal map $\Delta_{m+1} : X \to \prod^{m+1}X$ can be compressed into the `fat wedge' $\fatvee^{m+1}(X)$ (see Chapter X of \cite{Whitehead:elements}).
Since $\prod^{m+1}X/\fatvee^{m+1}(X)$ is the $(m+1)$-fold smash product $\wedge^{m+1}X$, we have a weaker invariant $\wcat{X}$, {\it the weak L-S category of $X$}, given by the least number $m$ such that the reduced diagonal map $\widebar{\Delta}_{m+1} : X \to \wedge^{m+1}X$ is trivial.
Hence $\wcat{X} \leq \cat{X}$.

T. Ganea has also introduced a stronger invariant $\Cat{X}$, {\it the strong L-S category of $X$}, by the least number $m$ such that there is a covering of $X$ by $(m+1)$ open subsets each of which is contractible in itself.
Thus $\wcat{X} \leq \cat{X} \leq \Cat{X}$.

The weak and strong L-S categories usually give nice estimates of L-S category especially for manifolds.
Actually, we do not know any example of a closed manifold whose strong L-S, L-S and weak L-S categories are not the same.
The following problems are posed by Ganea \cite{Ganea:conjecture}:
\begin{itemize}
\item[i)](Problem 1)~
Determine the L-S category of a manifold.
\item[ii)](Problem 4)~
Describe the L-S category of a sphere-bundle over a sphere in terms of homotopy invariants of the characteristic map of the bundle.
\end{itemize}

Problem 1 has been studied by many authors, such as Singhof \cite{Singhof:ls-lie-group,Singhof:ls-lie-group-ii,Singhof:minimal-cover}, Montejano \cite{Montejano:singhof}, Schweizer \cite{Schweizer:ls-sp2}, Gomez-Larra\~naga and Gonzalez-Acu\~na \cite{GG:ls-cat_3m}, James and Singhof \cite{JS:ls-so5} and Rudyak \cite{Rudyak:ls-cat_mfds,Rudyak:ls-cat_mfds2}.
In particular for compact simply-connected simple Lie groups, $\cat(SU(n+1))=n$ for $n \geq 1$ by \cite{Singhof:ls-lie-group}, $\cat(Sp(2))=3$ by \cite{Schweizer:ls-sp2} and $\cat(Sp(n)) \geq n+1$ for $n \geq 2$ by \cite{Singhof:ls-lie-group-ii}.
It was also announced recently that Problem 4 was solved by the first author \cite{Iwase:ls-sphere-bundle-over-sphere}.

%
The method in the present paper also provides a result for $G_2$, and thus we have the following result.
\begin{Thm}\label{thm:table}~
The following is the complete list of L-S categories of a simply-connected compact simple Lie group of rank $\leq 2$:
\begin{center}
\begin{tabular}{|c|c|c|c|c|}
\hline
{\small Lie groups} & 
{\small $Sp(1)=SU(2)=Spin(3)$}&
{\small $SU(3)$}&
{\small $Sp(2)=Spin(5)$}&
{\small $G_2$}
\\\hline
$\wcat{}$ & 
$1$&
$2$&
$3$&
$4$
\\\hline
$\cat{}$ & 
$1$&
$2$&
$3$&
$4$
\\\hline
$\Cat{}$ & 
$1$&
$2$&
$3$&
$4$
\\\hline
\end{tabular}
\end{center}
\end{Thm}
\par
Although the above result is known for experts, we give a short proof for $G_2$.
In fact, the result for $G_2$ has never been published and is obtained in a similar but easier manner than the following result for $Sp(3)$:
\begin{Thm}\label{thm:main}~
$\wcat(Sp(3)) = \cat(Sp(3)) = \Cat(Sp(3)) = 5$.
\end{Thm}
\begin{Rem}
The argument given to prove Theorem \ref{thm:main} provides an alternative proof of Schweizer's result $$\wcat(Sp(2)) = \cat(Sp(2)) = \Cat(Sp(2)) = 3.$$
\end{Rem}
The authors know that a similar result to Theorem \ref{thm:main} is obtained by Luc\'ia Fern\'andez-Su\'arez, Antonio G\'omez-Tato, Jeffrey Strom and Daniel Tanr\'e \cite{FGST:ls-cat_sp3}.
Our method is, however, much simpler and providing the following general result:
\begin{Thm}\label{thm:sp(n)}~
$n+2 \leq \wcat(Sp(n)) \leq \cat(Sp(n)) \leq \Cat(Sp(n))$ for $n \geq 3$.
\end{Thm}
This improves Singhof's result: $\cat(Sp(n)) \geq n+1$ for $n \geq 2$.
We propose the following conjecture.
\begin{Conj}~
Let $G$ be a simply-connected compact Lie group with $G$ = $\prod^{n}_{i=1} H_i$ where $H_i$ is a simple Lie group.
Then $\wcat(G)$ = $\cat(G)$ = $\Cat(G)$ and $\cat(G)$ = $\sum^{n}_{i=1} \cat(H_i)$.
\end{Conj}
It might be difficult to say something about $\cat{Sp(n)}$, but an old conjecture says the following.
\begin{Conj}~
$\cat{Sp(n)} = 2n-1$ for all $n \geq 1$.
\end{Conj}

The authors thank John Harper for many helpful conversations.

\section{Proof of Theorem \ref{thm:table}}
\par\noindent
%
Let us recall a CW decomposition of $G_2$ from \cite{MN:CW-decomp-G2}:
\begin{equation*}
G_2 = e^{0} \cup e^{3} \cup e^{5} \cup e^{6} \cup e^{8} \cup e^{9} \cup e^{11} \cup e^{14}.
\end{equation*}
On the other hand, we have the following cone-decomposition.
\begin{Thm}\label{thm:cone-decomp-G2}~
There is a cone-decomposition of $G_2$ as follows:
\begin{align*}&
G_2^{(5)} = {\Sigma}{\complex}P^2,\quad
S^{5} \cup e^{7} \to G_2^{(5)} \hookrightarrow G_2^{(8)},\quad
\\&
S^{8} \cup e^{10} \to G_2^{(8)} \hookrightarrow G_2^{(11)},\quad
S^{13} \to G_2^{(11)} \hookrightarrow G_2.\quad
\end{align*}
\end{Thm}
\begin{Proof}~
The first and the last formulae are obvious.
So we show the 2nd and 3rd formulae:
By taking the homotopy fibre $F_1$ of $G_2^{(5)} \hookrightarrow G_2$, we can easily observe using the Serre spectral sequence that the fibre has a CW structure given by $S^{5} \cup e^{7} \cup \mbox{(cells in dimensions $\geq 7$)}$, where the cohomology generators corresponding to $S^{5}$ and $e^{7}$ are transgressive.
Thus the mapping cone of $S^{5} \cup e^{7} \subset F_1 \to G_2^{(5)}$ has the homotopy type of $G_2^{(8)}$.
Similarly, the homotopy fibre $F_2$ of $G_2^{(8)} \hookrightarrow G_2$ has a CW structure given by $S^{8} \cup e^{10} \cup \mbox{(cells in dimensions $\geq 10$)}$, where the cohomology generators corresponding to $S^{8}$ and $e^{10}$ are transgressive.
Thus the mapping cone of $S^{8} \cup e^{10} \subset F_2 \to G_2^{(8)}$ has the homotopy type of $G_2^{(11)}$.
\end{Proof}
\begin{Cor}\label{cor:cone-decomp-G2}~
$1$ $\geq$ $\Cat(G_2^{(5)})$ $\geq$ $\Cat(G_2^{(3)})$, $2$ $\geq$ $\Cat(G_2^{(8)})$ $\geq$ $\Cat(G_2^{(6)})$, $3$ $\geq$ $\Cat(G_2^{(11)})$ $\geq$ $\Cat(G_2^{(9)})$ and $4$ $\geq$ $\Cat(G_2)$.
\end{Cor}
Let us recall the following well-known fact due to Borel.
\begin{Fact}\label{fact:cohomology-G2}~
$H^*(G_2;\integral/{2}\integral) \cong \integral/{2}\integral[x_3,x_5]/(x^4_3,x^2_5)$.
\end{Fact}
\begin{Cor}\label{cor:cohomology-G2}~
$\wcat(G_2^{(5)})$ $\geq$ $\wcat(G_2^{(3)})$ $\geq$ $1$, $\wcat(G_2^{(8)})$ $\geq$ $\wcat(G_2^{(6)})$ $\geq$ $2$, $\wcat(G_2^{(11)})$ $\geq$ $\wcat(G_2^{(9)})$ $\geq$ $3$ and $\wcat(G_2)$ $\geq$ $4$.
\end{Cor}
Corollaries \ref{cor:cone-decomp-G2} and \ref{cor:cohomology-G2} yield the following.
\begin{Thm}~
\begin{center}
\begin{tabular}{|c|c|c|c|c|c|c|c|}
\hline
{\small Skeleta}&%
{\small $G_2^{(3)}$}&%
{\small $G_2^{(5)}$}&%
{\small $G_2^{(6)}$}&%
{\small $G_2^{(8)}$}&%
{\small $G_2^{(9)}$}&%
{\small $G_2^{(11)}$}&%
{\small $G_2$}
\\\hline
$\wcat{}$ & 
$1$&
$1$&
$2$&
$2$&
$3$&
$3$&
$4$
\\\hline
$\cat{}$ & 
$1$&
$1$&
$2$&
$2$&
$3$&
$3$&
$4$
\\\hline
$\Cat{}$ &
$1$&
$1$&
$2$&
$2$&
$3$&
$3$&
$4$
\\\hline
\end{tabular}
\end{center}
\end{Thm}
\par\noindent
This completes the proof of Theorem \ref{thm:table}.
\section{The ring structure of $h^*(Sp(3))$}
\par\noindent
To show Theorem \ref{thm:main}, we introduce a cohomology theory $h^*(-)$ such that $h^*(X,A)$ $=$ $\{X/A,{\mathbb S}[0,2]\}$, where ${\mathbb S}[0,2]$ is the spectrum obtained from ${\mathbb S}$ by killing all homotopy groups of dimensions bigger than $2$.
Then ${\mathbb S}[0,2]$ is a ring spectrum with $\pi^S_*({\mathbb S}[0,2]) \cong \integral[\eta]/(\eta^{3},2\eta)$, where $\eta$ is the Hopf element in $\pi^S_{1}({\mathbb S}) = \pi^S_{1}({\mathbb S}[0,2])$.
Thus $h^*$ is an additive and multiplicative cohomology theory with $h^* = h^*(pt) \cong \integral[\varepsilon]/(\varepsilon^{3},2\varepsilon)$, $\deg\varepsilon=-1$, where $\varepsilon \in h^{-1} = \pi^S_{0}({\Sigma}^{-1}{\mathbb S}) \cong \pi^S_{1}({\mathbb S})$ corresponds to $\eta$.

The characteristic map of the principal $Sp(1)$-bundle 
\begin{equation*}
Sp(1) \hookrightarrow Sp(2) \rightarrow S^{7}
\end{equation*}
is given by $\omega=\langle\iota_3,\iota_3\rangle : S^{6} \to Sp(1) \homeo S^{3}$ the Samelson product of two copies of the identity $\iota_3 : S^3 \to S^3$, which is a generator of $\pi_{6}(S^{3}) \cong \integral/12\integral$.
We state the following well-known fact (see Whitehead \cite{Whitehead:elements}).
\begin{Fact}\label{prop:nu'}~
Let $\mu : S^3{\times}S^3 \to S^3$ be the multiplication of $Sp(1) \homeo S^3$.
Then we have 
\begin{align*}&
Sp(2) \simeq S^3 \cup_{\mu{\comp}(1{\times}\omega)} S^3{\times}C(S^{6})
= S^3 \cup_{\omega} C(S^{6}) \cup_{\hat{\mu}{\comp}[\iota_3,\omega]^r} C(S^{9}),
\end{align*}
where $\hat{\mu} : S^3{\times}S^3 \cup_{{\ast}{\times}{\omega}} \{\ast\}{\times}C(S^{6}) \to S^{3} \cup_{\omega} C(S^{6})$ is given by $\hat{\mu}\vert_{S^3{\times}S^3} = \mu$ and $\hat{\mu}\vert_{S^{3} \cup_{\omega} C(S^{6})} = 1$ the identity and $[\iota_3,\chi_{\omega}]^r : S^{9} \to S^{3}{\times}S^3 \cup_{{\ast}{\times}{\omega}} \{\ast\}{\times}C(S^{6})$ is the relative Whitehead product of the identity $\iota_3 : S^{3} \to S^{3}$ and the characteristic map $\chi_{\omega} : (C(S^{6}),S^{6}) \to (S^{3} \cup e^{7},S^{3})$ of the $7$-cell.
Thus we have $1 \geq \Cat(Sp(2)^{(3)})$, $2 \geq \Cat(Sp(2)^{(7)})$ and $3 \geq \Cat(Sp(2))$.
\end{Fact}
Let $\nu : S^{7} \to S^{4}$ be the Hopf element whose suspension $\nu_{n}={\Sigma}^{n-4}\nu$ ($n \geq 4$) gives a generator of $\pi_{n+3}(S^{n}) \cong \integral/24\integral$ for $n \geq 5$.
Then we remark that $\omega_{n}={\Sigma}^{n-3}\omega$ ($n \geq 3$) satisfies the formula $\omega_{n}=2\nu_{n} \in \pi_{n+3}(S^{n})$ for $n \geq 5$.
By Zabrodsky \cite{Zabrodsky:Hopf-space}, there is a natural splitting 
$${\Sigma}(S^{3}{\times}S^{3} \cup \{\ast\}{\times}(S^{3} \cup_{\omega} e^{7})) \simeq {\Sigma}S^{3} \vee {\Sigma}(S^{3} \cup_{\omega} e^{7}) \vee {\Sigma}S^{3}{\wedge}S^{3}.$$
Then by the definition of a relative Whitehead product, the composition of $[\iota_3,\omega]^r$ with the projections to $S^{3}$ and $S^{3} \cup_{\omega} e^{7}$ are trivial and the composition with the projection to $S^{3}{\wedge}S^{3}$ is given by $\iota_3{\wedge}\omega$.
Thus we have 
$${\Sigma}(\hat{\mu}{\comp}[\iota_3,\omega]^r) = H(\mu){\comp}{\Sigma}(\iota_3{\wedge}\omega) = \pm \nu{\comp}\omega_{7} = 2\nu{\comp}\nu_{7} \not= 0$$
in $\pi_{10}(S^{4}) \cong \integral/24\integral\langle\nu{\comp}\nu_{7}\rangle{\oplus}\integral/2\integral\langle\omega_{4}{\comp}\nu_{7}\rangle$, and hence we have 
$${\Sigma}^2(\hat{\mu}{\comp}[\iota_3,\omega]^r) = \nu_5{\comp}\omega_8 = 2\nu^2_5 = 0 \in \pi_{11}(S^{5}) \cong \integral/2\integral$$
by Proposition 5.11 of Toda \cite{Toda:comp-methods}.
The following two facts are also well-known.
\begin{Fact}\label{fact:Sigma2Sp(2)}~
We have the following homotopy equivalences:
\begin{align*}&
Sp(2)/S^3 \simeq (S^3{\times}C(S^{6}))/(S^3{\times}S^{6})
= S^3_{+}{\wedge}{\Sigma}(S^{6})
= S^{7} \vee S^{10},
\\&
{\Sigma}^2Sp(2) \simeq {\Sigma}^2(S^{3} \cup_{\omega} C(S^{6})) \vee {\Sigma}^2S^{10} = S^{5} \cup_{\omega_{5}} C(S^{8}) \vee S^{12}.
\end{align*}
\end{Fact}
\par
\begin{Fact}\label{prop:nu}~
The $11$-skeleton $X_{3,2}^{(11)}$ of $X_{3,2} = Sp(3)/Sp(1)$ has the homotopy type of $S^{7} \cup_{\nu_{7}} e^{11}$.
\end{Fact}
Restricting the principal $Sp(1)$-bundle %
$
Sp(1) \hookrightarrow Sp(3) \overset{q}\rightarrow X_{3,2}
$ 
to the subspace $X_{3,2}^{(11)} = S^{7} \cup_{\nu_{7}} e^{11}$ of $X_{3,2}$, we obtain the subspace $q^{-1}(X_{3,2}^{(11)}) = Sp(3)^{(14)}$ of $Sp(3)$ as the total space of the principal $Sp(1)$-bundle %
$
Sp(1) \hookrightarrow Sp(3)^{(14)} \overset{q}\rightarrow {\Sigma}(S^{6} \cup_{\nu_{6}} e^{10})
$ 
with a characteristic map $\phi : S^{6} \cup_{\nu_{6}} e^{10} \to Sp(1) \homeo S^3$, which is an extension of $\omega : S^{6} \to S^{3}$.
\begin{Prop}\label{prop:phi}~
We have the following homotopy equivalences:
\begin{align*}&
Sp(3)^{(14)} \simeq S^3 \cup_{\mu{\comp}(1{\times}\phi)} S^3{\times}C(S^{6} \cup_{\nu_{6}} e^{10})
\\&\hspace{14.0mm}
= S^3 \cup_{\phi} C(S^{6} \cup_{\nu_{6}} e^{10}) \cup C(S^{9} \cup_{\nu_{9}} e^{13}),
\\&
Sp(3)^{(14)}/S^3 \simeq (S^3{\times}C(S^{6} \cup_{\nu_{6}} e^{10}))/(S^3{\times}(S^{6} \cup_{\nu_{6}} e^{10}))
\\&\hspace{19.65mm}
= S^3_{+}{\wedge}{\Sigma}(S^{6} \cup_{\nu_{6}} e^{10})
= (S^{7} \cup_{\nu_{7}} e^{11}) \vee (S^{10} \cup_{\nu_{10}} e^{14}),
\\&
Sp(n) \simeq Sp(n-1) \cup Sp(n-1){\times}C(S^{4n-2}),
\quad \intertext{where $Sp(n-1) \subset Sp(n)^{((2n+1)n-11)}$ for $n \geq 3$, and hence}
&
Sp(n)/Sp(n)^{((2n+1)n-11)} 
\\&\qquad
\simeq (Sp(n-1){\times}C(S^{4n-2}))/(Sp(n-1){\times}S^{4n-2} 
\\&\qquad\qquad\qquad\qquad\qquad\qquad\qquad\ 
\cup Sp(n-1)^{((2n-1)(n-1)-11)}{\times}C(S^{4n-2}))
\\&\qquad
= (Sp(n-1)/Sp(n-1)^{((2n-1)(n-1)-11)}){\wedge}{\Sigma}S^{4n-2}
\\&\qquad
= \cdots = (Sp(2)/\emptyset){\wedge}{\Sigma}S^{10}{\wedge}\cdots{\wedge}{\Sigma}S^{4n-2}
= (Sp(2)_{+}){\wedge}S^{(2n+1)n-10}
\\&\qquad
= S^{(2n+1)n-10} \vee S^{(2n+1)n-10}{\wedge}Sp(2) 
\\&\qquad
= S^{(2n+1)n-10} \vee (S^{(2n+1)n-7} \cup_{\omega_{(2n+1)n-7}} e^{(2n+1)n-3}) \vee S^{(2n+1)n},\quad\text{for $n \geq 3$.}
\end{align*}
\end{Prop}
This yields the following result.
\begin{Prop}\label{prop:cw-decomp}~
Let $\hat{\mu} : S^3{\times}S^3 \cup_{{\ast}{\times}{\phi}} \{\ast\}{\times}(S^{3} \cup_{\phi} C(S^{6} \cup_{\nu_6} e^{10})) \to S^{3} \cup_{\phi} C(S^{6} \cup_{\nu_6} e^{10})$ be the map given by $\hat{\mu}\vert_{S^3{\times}S^3} = \mu$ 
and $\hat{\mu}\vert_{S^{3} \cup_{\phi} C(S^{6} \cup_{\nu_6} e^{10})} = 1$ the identity.
Then we have the following cone decomposition of $Sp(3)$:
\begin{equation*}
Sp(3) \simeq S^{3} \cup_{\phi} C(S^{6} \cup_{\nu_6} e^{10}) \cup_{\hat{\mu}{\comp}\hat{\phi}} C(S^{9} \cup_{\nu_9} e^{13}) \cup C(S^{17}) \cup C(S^{20}).
\end{equation*}
\end{Prop}
\begin{Cor}\label{cor:cw-decomp}~
$1$ $\geq$ $\Cat(Sp(3)^{(3)})$, $2$ $\geq$ $\Cat(Sp(3)^{(7)})$, $3$ $\geq$ $\Cat(Sp(3)^{(14)})$ $\geq$ $\Cat(Sp(3)^{(11)})$ $\geq$ $\Cat(Sp(3)^{(10)})$, $4$ $\geq$ $\Cat(Sp(3)^{(18)})$ and $5$ $\geq$ $\Cat(Sp(3))$.
\end{Cor}
To determine the ring structures of $h^*(Sp(2))$ and $h^*(Sp(3))$, we show the following lemma.
\begin{Lem}\label{thm:h-inv}
Let $h^*$ be any multiplicative generalised cohomology theory and let $Q = S^{r} \cup_{f} e^{q}$ for a given map $f : S^{q-1} \to S^{r}$ with $h^*(Q) \cong h^*\langle{1,x,y}\rangle$, where $x$ and $y$ correspond to the generators of $h^*(S^{r}) \cong h^*\langle{x_0}\rangle$ and $h^*(S^{q}) \cong h^*\langle{y_0}\rangle$.
Then 
$$x^2 = \pm H^h(\lambda_2(f)){\cdot}y \quad\text{in}\quad h^{\ast}(Q),$$
where $\lambda_2(f) \in \pi_{q}(S^{2r})$ is the Boardman-Steer Hopf invariant equal to ${\Sigma}h_2(f)$ the suspension of the James-Hopf invariant $h_2(f)$ (see \cite{BS:hopf-invariants}) and $H^h : \pi_{q}(S^{2r}) \to h^{2r}(S^{q}) \cong h^{2r-q}$ is the Hurewicz homomorphism given by $H^h(g) = \Sigma_{\ast}^{-q}g^{\ast}(x_0{\otimes}x_0)$.
\end{Lem}
\begin{Proof}~
By Boardman and Steer \cite{BS:hopf-invariants}, $\widebar{\Delta} : Q_2=S^{r} \cup_{f} e^{q} \to Q_2{\wedge}Q_2$ equals the composition $(i_2{\wedge}i_2){\comp}\lambda_2(f){\comp}q_2$, where $q_2 : Q_2 \to Q_2/S^{r}=S^{q}$ is the collapsing map and $i_2 : S^{r} \hookrightarrow Q_2$ is the bottom-cell inclusion.
Thus we have 
\begin{align*}
x^2 
&= \widebar{\Delta}^*(x{\otimes}x) 
= ((i_2{\wedge}i_2){\comp}\lambda_2(f){\comp}q_2)^*(x{\otimes}x) 
\\&
= q_2^*(\lambda_2(f)^*(i_2^*(x){\otimes}i_2^*(x))) 
= q_2^*(\lambda_2(f)^*(x_0{\otimes}x_0)) 
= q_2^*(\Sigma_{\ast}^{q}H^h(\lambda_2(f))).
\intertext{Since $\Sigma_{\ast}^{q}H^h(\lambda_2(f))$ is $H^h(\lambda_2(f)){\cdot}y_0 \in h^{2r}(S^q)$ up to sign, we proceed as} x^2 &
= q_2^*(\pm H^h(\lambda_2(f)){\cdot}y_0) 
= \pm H^h(\lambda_2(f)){\cdot}q_2^*(y_0) 
= \pm H^h(\lambda_2(f)){\cdot}y.
\end{align*}
This completes the proof of the lemma.
\end{Proof}

Using cohomology long exact sequences derived from the cell structure of $Sp(3)$ and a direct calculation using Proposition \ref{prop:phi} and Lemma \ref{thm:h-inv}, we deduce the following result for the cohomology theory $h^*$ considered at the beginning of this section.
\begin{Thm}\label{prop:main}~
The ring structures of $h^*(Sp(2))$ and $h^*(Sp(3))$ are as follows:
\begin{align*}&
h^*(Sp(2)) \cong h^*\{1,x_3,x_7,y_{10}\},
\\&
h^*(Sp(3)) \cong h^*\{1,x_3,x_7,x_{11},y_{10},y_{14},y_{18},z_{21}\}
\end{align*}
%
with ring structures given by $x_{3}^2$ = $\varepsilon{\cdot}x_{7}$, $x_{7}^2$ = $0$, $x_{11}^2$ = $0$, $x_{3}x_{7}$ = $y_{10}$, $x_{3}x_{11}$ = $y_{14}$, $x_{7}x_{11}$ = $y_{18}$ and $x_{3}x_{7}x_{11}$ = $z_{21}$, where $\varepsilon$ is the non-zero element in $h^{-1}$.
\end{Thm}
\begin{Rem}\label{rem:attaching_map}~
The two possible attaching maps $: S^{10} \to S^3 \cup_{\omega} e^7$ of $e^{11}$ discovered by Luc\'ia Fern\'andez-Su\'arez, Antonio G\'omez-Tato and Daniel Tanr\'e \cite{FGT:ls-cat_sp3} are homotopic in $Sp(2)$.
So, we can not find any effective difference in the ring structure of $h^*(Sp(3))$ by altering, as is performed in \cite{FGST:ls-cat_sp3}, the attaching map of $e^{11}$.
\end{Rem}
\begin{Cor}~
$\wcat(Sp(3)^{(3)})$ $\geq$ $1$, $\wcat(Sp(3)^{(7)})$ $\geq$ $2$, $\wcat(Sp(3)^{(18)})$ $\geq$ $\wcat(Sp(3)^{(14)})$ $\geq$ $\wcat(Sp(3)^{(11)})$ $\geq$ $\wcat(Sp(3)^{(10)})$ $\geq$ $3$ and $\wcat(Sp(3))$ $\geq$ $4$, together with $\wcat(Sp(2)^{(3)})$ $\geq$ $1$, $\wcat(Sp(2)^{(7)})$ $\geq$ $2$ and $\wcat(Sp(2))$ $\geq$ $3$.
\end{Cor}
\begin{Cor}~
\begin{center}
\begin{tabular}{|c|c|c|c|}
\hline
{\small Skeleta}&%
{\small $Sp(2)^{(3)}$}&%
{\small $Sp(2)^{(7)}$}&%
{\small $Sp(2)$}
\\\hline
$\wcat{}$ & 
$1$&
$2$&
$3$
\\\hline
$\cat{}$ & 
$1$&
$2$&
$3$
\\\hline
$\Cat{}$ &
$1$&
$2$&
$3$
\\\hline
\end{tabular}
\end{center}
\end{Cor}
%
\section{Proof of Theorem \ref{thm:main}}
\par\noindent
By Facts \ref{prop:nu'} and \ref{fact:Sigma2Sp(2)}, %
the smash products ${\wedge}^{4}Sp(3)$ and ${\wedge}^{5}Sp(3)$ satisfy
\begin{align*}&
({\wedge}^{4}Sp(3))^{(19)} \simeq S^{12} \cup_{\omega_{12}} e^{16} \vee (S^{16} \vee S^{16} \vee S^{16}) \vee (S^{19} \vee S^{19} \vee S^{19} \vee S^{19}),
\\&
({\wedge}^{5}Sp(3))^{(22)} \simeq S^{15} \cup_{\omega_{15}} e^{19} \vee (S^{19} \vee S^{19} \vee S^{19}) \vee (S^{22} \vee S^{22} \vee S^{22} \vee S^{22}).
\end{align*}
Then we have the following two propositions.
\begin{Prop}\label{prop:wedge-sp(3)}~
The bottom-cell inclusions $i : S^{12} \hookrightarrow {\wedge}^{4}Sp(3)$ and $i' : S^{15} \hookrightarrow {\wedge}^{5}Sp(3)$ induce injective homomorphisms
\begin{align*}&
i_{\ast} : \pi_{18}(S^{12}) \to \pi_{18}({\wedge}^{4}Sp(3)^{(18)})\quad\mbox{and}\quad
i'_{\ast} : \pi_{21}(S^{15}) \to \pi_{21}({\wedge}^{5}Sp(3)),
\end{align*}
respectively.
\end{Prop}
\begin{Proof}~
We have the following two exact sequences
\begin{align*}&
\pi_{18}(S^{15}) \overset{\psi}\to \pi_{18}(S^{12}) \overset{i_{\ast}}\to \pi_{18}({\wedge}^{4}Sp(3)^{(18)}) \to \pi_{18}(S^{16}{\vee}S^{16}{\vee}S^{16}{\vee}S^{16}),
\\&
\pi_{21}(S^{18}) \overset{\psi'}\to \pi_{21}(S^{15}) \overset{i'_{\ast}}\to \pi_{21}({\wedge}^{5}Sp(3)) \to \pi_{21}(S^{19}{\vee}S^{19}{\vee}S^{19}{\vee}S^{19}{\vee}S^{19}),
\end{align*}
where $\pi_{18}(S^{12}) \cong \pi_{21}(S^{15}) \cong \integral/2\integral\nu^2_{15}$ and $\psi$ and $\psi'$ are induced from $\omega_{12}=2\nu_{12}$ and $\omega_{15}=2\nu_{15}$.
Thus $\psi$ and $\psi'$ are trivial, and hence $i_{\ast}$ and $i'_{\ast}$ are injective.
\end{Proof}
\begin{Prop}\label{prop:injective}~
The collapsing maps $q : Sp(3)^{(18)} \to Sp(3)^{(18)}/Sp(3)^{(14)}$ = $S^{18}$ and $q' : Sp(3) \to Sp(3)/Sp(3)^{(18)}$ = $S^{21}$ induce injective homomorphisms
\begin{enumerate}
\item[]
$q^{\ast} : \pi_{18}({\wedge}^4Sp(3)^{(18)}) \to [Sp(3)^{(18)},{\wedge}^4Sp(3)^{(18)}]$\quad{and}\quad 
\item[]
${q'}^{\ast} : \pi_{21}({\wedge}^5Sp(3)) \to [Sp(3),{\wedge}^5Sp(3)]$, 
\end{enumerate}
respectively.
\end{Prop}
\begin{Proof}~
Firstly, we show that ${q'}^{\ast}$ is injective:
Since we have $[Sp(3),{\wedge}^5Sp(3)]$ = $[(S^{14} \cup_{\omega_{14}} e^{18}) \vee S^{21},{\wedge}^5Sp(3)]$ = $[S^{14} \cup_{\omega_{14}} e^{18},{\wedge}^5Sp(3)]{\oplus}\pi_{21}({\wedge}^5Sp(3))$ by Proposition \ref{prop:phi}, ${q'}^{\ast}$ is clearly injective.
\par
Secondly, we show that ${q}^{\ast}$ is injective:
Similarly we have $[Sp(3)^{(18)},{\wedge}^4Sp(3)^{(18)}]$ = $[S^{14} \cup_{\omega_{14}} e^{18},{\wedge}^4Sp(3)^{(18)}]$ by Proposition \ref{prop:phi}.
Thus it is sufficient to show that $\bar{q}^{\ast} : \pi_{18}({\wedge}^4Sp(3)^{(18)}) \to [S^{14} \cup_{\omega_{14}} e^{18},{\wedge}^4Sp(3)^{(18)}]$ is injective, where $\bar{q} : S^{14} \cup_{\omega_{14}} e^{18} \to S^{18}$ is the collapsing map.
In the exact sequence 
\begin{equation*}
\pi_{15}({\wedge}^4Sp(3)^{(18)}) \overset{{\omega_{15}}^{\ast}}\to \pi_{18}({\wedge}^4Sp(3)^{(18)}) \overset{\bar{q}^{\ast}}\to [S^{14} \cup_{\omega_{14}} e^{18},{\wedge}^4Sp(3)^{(18)}],
\end{equation*}
we know that $\pi_{15}({\wedge}^4Sp(3)^{(18)}) \cong \pi_{15}(S^{12} \cup_{\omega_{12}} e^{16}) = \integral/2\integral$ is generated by the composition of $\nu_{12}$ and the bottom-cell inclusion.
Since $\nu_{12}{\comp}\omega_{15} = 0 \in \pi_{18}(S^{12})$, the homomorphism ${\omega_{15}}^{\ast}$ is trivial, and hence $\bar{q}^{\ast}$ is injective.
\end{Proof}
Then the following lemma implies that $\widebar{\Delta}_{4}$ and $\widebar{\Delta}_{5}$ are non-trivial by Propositions \ref{prop:wedge-sp(3)} and \ref{prop:injective}.
\begin{Lem}\label{lem:key}~
We obtain that $\widebar{\Delta}_{4}$ = $i{\comp}\nu_{12}^2{\comp}q : Sp(3)^{(18)} \to {\wedge}^{4}Sp(3)^{(18)}$ and that $\widebar{\Delta}_{5}$ = $i'{\comp}\nu_{15}^2{\comp}q' : Sp(3) \to {\wedge}^{5}Sp(3)$.
\end{Lem}
\begin{Proof}~
Firstly, we show that $\widebar{\Delta}_{4} = i{\comp}\nu_{12}^2{\comp}q$ implies $\widebar{\Delta}_{5} = i'{\comp}\nu_{15}^2{\comp}q'$.
For dimensional reasons, the image of $\widebar{\Delta} : Sp(3) \to Sp(3){\wedge}Sp(3)$ is in $Sp(3)^{(18)}{\wedge}Sp(3)^{(14)} \cup S^{3}{\wedge}Sp(3)^{(18)}$.
Since $Sp(3)^{(14)}$ is of cone-length $3$ by Corollary \ref{cor:cw-decomp}, the restriction of the map $(1{\wedge}\widebar{\Delta}_4){\comp}\widebar{\Delta}=\widebar{\Delta}_5$ to $Sp(3)^{(18)}{\wedge}Sp(3)^{(14)}$ is trivial.
Thus $\widebar{\Delta}_5$ equals the composition
\begin{equation*}
\widebar{\Delta}_5 : Sp(3) \to S^{3}{\wedge}Sp(3)^{(18)} \overset{1{\wedge}\widebar{\Delta}_{4}}\to {\wedge}^{5}Sp(3)^{(18)} \subset {\wedge}^{5}Sp(3).
\end{equation*}
Then by $\widebar{\Delta}_{4} = i{\comp}\nu_{12}^2{\comp}q$, we observe that $\widebar{\Delta}_5 = i'{\comp}(\iota_3{\wedge}\nu_{12}^2){\comp}q' = i'{\comp}\nu_{15}^2{\comp}q'$.
\par
So, we are left to show $\widebar{\Delta}_{4} = i{\comp}\nu_{12}^2{\comp}q$.
For dimensional reasons, the image of $\widebar{\Delta} : Sp(3)^{(18)} \to Sp(3)^{(18)}{\wedge}Sp(3)^{(18)}$ is in $Sp(3)^{(14)}{\wedge}S^{3} \cup Sp(3)^{(11)}{\wedge}Sp(3)^{(7)} \cup Sp(3)^{(7)}{\wedge}Sp(3)^{(11)} \cup S^{3}{\wedge}Sp(3)^{(14)}$.
Since $S^{3} \cup_{\phi} C(S^{6} \cup_{\nu_{6}} e^{10})$ is of cone-length $2$ by Corollary \ref{cor:cw-decomp}, the restriction of $\widebar{\Delta}_{3} : Sp(3)^{(18)} \to {\wedge}^{3}Sp(3)^{(18)}$ to $S^{3} \cup_{\phi} C(S^{6} \cup_{\nu_{6}} e^{10})$ is trivial.
Hence $1{\wedge}\widebar{\Delta}_{3} : Sp(3)^{(14)}{\wedge}S^{3} \cup Sp(3)^{(11)}{\wedge}Sp(3)^{(7)} \cup Sp(3)^{(7)}{\wedge}Sp(3)^{(11)} \cup S^{3}{\wedge}Sp(3)^{(14)} \to {\wedge}^{4}Sp(3)^{(18)}$ equals the composition
\begin{align*}
1{\wedge}\widebar{\Delta}_{3} : (Sp(3){\wedge}Sp(3))^{(18)}
&\overset{\alpha}\to (S^{3} \cup_{\omega} e^{7}){\wedge}S^{10} \cup
S^{3}{\wedge}(S^{10} \cup_{\nu_{10}} e^{14}) 
\overset{1{\wedge}\beta}\to {\wedge}^{4}(S^{3} \cup_{\omega} e^{7}).
\end{align*}
The map $\alpha{\comp}\widebar{\Delta} : Sp(3)^{(18)} \to (S^{3} \cup_{\omega} e^{7}){\wedge}S^{10} \cup S^{3}{\wedge}(S^{10} \cup_{\nu_{10}} e^{14})$ equals the composition
\begin{equation*}
\alpha{\comp}\widebar{\Delta} : Sp(3)^{(18)} \to S^{14} \cup_{\omega_{14}} e^{18} \to (S^{3} \cup_{\omega} e^{7}){\wedge}S^{10} \cup S^{3}{\wedge}(S^{10} \cup_{\nu_{10}} e^{14}).
\end{equation*}
Collapsing the subspace $S^{3}{\wedge}(S^{10} \cup_{\nu_{10}} e^{14})$ of $(S^{3} \cup_{\omega} e^{7}){\wedge}S^{10} \cup S^{3}{\wedge}(S^{10} \cup_{\nu_{10}} e^{14})$, we obtain a map 
$$
q'{\comp}\alpha{\comp}\widebar{\Delta} : Sp(3)^{(18)} \to S^{7}{\wedge}S^{10},
$$
where $q' : (S^{3} \cup_{\omega} e^{7}){\wedge}S^{10} \cup S^{3}{\wedge}(S^{10} \cup_{\nu_{10}} e^{14}) \to S^{3}{\wedge}S^{10}$ is the collapsing map.
For dimensional reasons, $q'{\comp}\alpha{\comp}\widebar{\Delta}$ equals the composition:
\begin{equation*}
q'{\comp}\alpha{\comp}\widebar{\Delta} : Sp(3)^{(18)} \to Sp(3)^{(18)}/Sp(3)^{(14)}=S^{18} \overset{\gamma}\to S^{7}{\wedge}S^{10}.
\end{equation*}
If $\gamma$ were non-trivial, then $\gamma$ would be $\eta_{17} : S^{18} \to S^{17}$, and hence we should have $x_{7}y_{10} = \varepsilon{\cdot}y_{18} \not= 0$.
However, from the ring structure of $h^*(Sp(3))$ given in Theorem \ref{prop:main}, we know $x_{7}y_{10} = 0$, and hence we obtain $\gamma=0$.
Then the image of $\alpha{\comp}\widebar{\Delta}$ is in the subspace $S^{3}{\wedge}(S^{10} \cup_{\nu_{10}} e^{14})$ of $(S^{3} \cup_{\omega} e^{7}){\wedge}S^{10} \cup S^{3}{\wedge}(S^{10} \cup_{\nu_{10}} e^{14})$, since they are $12$-connected.
Hence $\widebar{\Delta}_{4} = (1{\wedge}\widebar{\Delta}_{3}){\comp}\widebar{\Delta}$ equals the composition
\begin{align*}&
\widebar{\Delta}_{4} :  Sp(3)^{(18)} \overset{\alpha{\comp}\widebar{\Delta}}\to S^{3}{\wedge}(S^{10} \cup_{\nu_{10}} e^{14}) 
\overset{1{\wedge}\beta}\to S^{3}{\wedge}({\wedge}^{3}(S^{3} \cup_{\omega} e^{7}))^{(15)} \subset {\wedge}^{4}Sp(3)^{(18)},
\end{align*}
where $({\wedge}^{3}(S^{3} \cup_{\omega} e^{7}))^{(15)}$ is given by 
$(S^{3} \cup_{\omega} e^{7}){\wedge}S^{3}{\wedge}S^{3} \,\cup\, S^{3}{\wedge}(S^{3} \cup_{\omega} e^{7}){\wedge}S^{3} \,\cup\, S^{3}{\wedge}S^{3}{\wedge}(S^{3} \cup_{\omega} e^{7}).$
Collapsing the subspace ${\wedge}^{3}S^{3}$ of $({\wedge}^{3}(S^{3} \cup_{\omega} e^{7}))^{(15)}$, we obtain a map 
$$
q''{\comp}\beta : S^{10} \cup_{\nu_{10}} e^{14} \to S^{7}{\wedge}S^{3}{\wedge}S^{3} \,\cup\, S^{3}{\wedge}S^{7}{\wedge}S^{3} \,\cup\, S^{3}{\wedge}S^{3}{\wedge}S^{7}, 
$$
where $q'' : ({\wedge}^{3}(S^{3} \cup_{\omega} e^{7}))^{(15)} \to S^{7}{\wedge}S^{3}{\wedge}S^{3} \,\cup\, S^{3}{\wedge}S^{7}{\wedge}S^{3} \,\cup\, S^{3}{\wedge}S^{3}{\wedge}S^{7}$ is the collapsing map.
For dimensional reasons, $q''{\comp}\beta$ equals the composition
\begin{equation*}
q''{\comp}\beta : S^{10} \cup_{\nu_{10}} e^{14} \to S^{14} \overset{\gamma'}\to S^{7}{\wedge}S^{3}{\wedge}S^{3} \,{\vee}\, S^{3}{\wedge}S^{7}{\wedge}S^{3} \,{\vee}\, S^{3}{\wedge}S^{3}{\wedge}S^{7}.
\end{equation*}
If $\gamma'$ were non-trivial, then its projection to $S^{13}$ would be $\eta_{13} : S^{14} \to S^{13}$, and hence we should have $x_{3}^2x_{7} = \varepsilon{\cdot}y_{14} \not= 0$.
However, from the ring structure of $h^*(Sp(3))$ given in Theorem \ref{prop:main}, we know $x_{3}^2x_{7} = \varepsilon{\cdot}x_{7}^2 = 0$, and hence we obtain $\gamma'=0$.
Hence the image of $\beta$ lies in the subspace ${\wedge}^{3}S^{3}$ of ${\wedge}^{3}Sp(3)^{(18)}$.
\par
On the other hand, for dimensional reasons, $\alpha{\comp}\widebar{\Delta}$ equals the composition
\begin{equation*}
\alpha{\comp}\widebar{\Delta} : Sp(3)^{(18)} \to S^{14} \cup_{\omega_{14}} e^{18} \overset{\alpha'}\to S^{3}{\wedge}(S^{10} \cup_{\nu_{10}} e^{14}),
\end{equation*}
where the restriction $\alpha'\vert_{S^{14}}$ equals the composition 
$$\alpha'\vert_{S^{14}} : S^{14} \overset{\gamma''}\to S^{13} \hookrightarrow S^{3}{\wedge}(S^{10} \cup_{\nu_{10}} e^{14}).$$
If it were non-trivial, then $\gamma''$ would be $\eta_{13} : S^{14} \to S^{13}$, and hence we should have $x_{3}y_{10} = \varepsilon{\cdot}y_{14} \not= 0$.
However, from the ring structure of $h^*(Sp(3))$ given in Theorem \ref{prop:main}, we know $x_{3}y_{10} = x_{3}^2x_{7} = \varepsilon{\cdot}x_{7}^2 = 0$, and hence we obtain $\gamma''=0$.
Hence $\alpha{\comp}\widebar{\Delta}$ equals the composition
\begin{equation*}
\alpha{\comp}\widebar{\Delta} : Sp(3)^{(18)} \overset{q}\to S^{18} \overset{\alpha''}\to S^{3}{\wedge}(S^{10} \cup_{\nu_{10}} e^{14}),
\end{equation*}
and hence $\widebar{\Delta}_{4}$ equals the composition
\begin{align*}
\widebar{\Delta}_{4} :  Sp(3)^{(18)} \overset{q}\to S^{18} &\overset{\alpha''}\to S^{3}{\wedge}(S^{10} \cup_{\nu_{10}} e^{14}) 
\overset{1{\wedge}\beta}\to S^{3}{\wedge}({\wedge}^{3}S^{3}) \overset{i}\hookrightarrow {\wedge}^{4}Sp(3)^{(18)}.
\end{align*}
Now, we are ready to determine $\widebar{\Delta}_{4}$:
By Theorem \ref{prop:main}, we know $x_{3}^2x_{11}=\varepsilon{\cdot}z_{18}$ and $x_3^2=\varepsilon{\cdot}x_7$, hence $\alpha'' : S^{18} \to S^{13} \cup_{\nu_{13}} e^{17}$ is a co-extension of $\eta_{16} : S^{17} \to S^{16}$ on $S^{13} \cup_{\nu_{13}} e^{17}$ and $1{\wedge}\beta : S^{13} \cup_{\nu_{13}} e^{17} \to S^{12}$ is an extension of $\eta_{12} : S^{13} \to S^{12}$.
Thus the composition $(1{\wedge}\beta){\comp}\alpha''$ is an element of the Toda bracket $\{\eta_{12},\nu_{13},\eta_{16}\}$ which contains a single element $\nu_{12}^{2}$ by Lemma 5.12 of  \cite{Toda:comp-methods}, and hence $\widebar{\Delta}_{4} = i{\comp}\nu_{12}^2{\comp}q$.
\end{Proof}
\begin{Cor}
$\wcat(Sp(3)^{(18)})$ $\geq$ $4$ and $\wcat(Sp(3))$ $\geq$ $5$.
\end{Cor}
This yields the following result.
\begin{Thm}~
\begin{center}
\begin{tabular}{|c|c|c|c|c|c|c|c|}
\hline
{\small Skeleta}&%
{\small $Sp(3)^{(3)}$}&%
{\small $Sp(3)^{(7)}$}&%
{\small $Sp(3)^{(10)}$}&%
{\small $Sp(3)^{(11)}$}&%
{\small $Sp(3)^{(14)}$}&%
{\small $Sp(3)^{(18)}$}&%
{\small $Sp(3)$}
\\\hline
$\wcat{}$ & 
$1$&
$2$&
$3$&
$3$&
$3$&
$4$&
$5$
\\\hline
$\cat{}$ & 
$1$&
$2$&
$3$&
$3$&
$3$&
$4$&
$5$
\\\hline
$\Cat{}$ &
$1$&
$2$&
$3$&
$3$&
$3$&
$4$&
$5$
\\\hline
\end{tabular}
\end{center}
\end{Thm}
\par\noindent
This completes the proof of Theorem \ref{thm:main}.

\section{Proof of Theorem \ref{thm:sp(n)}}
\par\noindent
We know that for $n \geq 4$,
\begin{align*}&
Sp(n)^{(16)} = Sp(4)^{(15)} = Sp(3)^{(14)} \cup e^{15},
\\&
Sp(n)^{(19)} = 
\left\{\begin{array}{ll}
Sp(4)^{(15)} \cup (e^{18} \vee e^{18})&n = 4,\\
Sp(4)^{(15)} \cup (e^{18} \vee e^{18}) \cup e^{19}&n \geq 5,
\end{array}\right.
\\&
Sp(n)^{(21)} = 
Sp(n)^{(19)} \cup e^{21}
\end{align*}
and that $\wcat(Sp(3)^{(14)}) = \cat(Sp(3)^{(14)}) = \Cat(Sp(3)^{(14)}) = 3$.
Firstly, we show the following.
\begin{Prop}\label{prop:Sp(n)-15}~
$\wcat(Sp(4)^{(15)}) = 3$.
\end{Prop}
\begin{Proof}~
Since $\cat(Sp(2))=3$, it follows that $\wcat(Sp(4)^{(15)}) \geq 3$ by Theorem 5.1 of \cite{Iwase:counter-ls-m}.
Hence we are left to show $\wcat(Sp(4)^{(15)}) \leq 3$:
For dimensional reasons, $\widebar{\Delta}_{4} = (\widebar{\Delta}{\wedge}\widebar{\Delta}){\comp}\widebar{\Delta} : Sp(4)^{(15)} \to {\wedge}^{4}Sp(4)^{(15)}$ equals the composition
\begin{equation*}
\widebar{\Delta}_{4} : Sp(4)^{(15)} \overset{\alpha_0}\to Sp(4)^{(11)}{\wedge}Sp(4)^{(11)} \overset{\widebar{\Delta}{\wedge}\widebar{\Delta}}\to {\wedge}^4Sp(4)^{(11)} \hookrightarrow {\wedge}^{4}Sp(4)^{(15)},
\end{equation*}
for some $\alpha_0$.
By Fact \ref{fact:Sigma2Sp(2)}, 
$\widebar{\Delta} : Sp(4)^{(11)} \to {\wedge}^2Sp(4)^{(11)}$ equals the composition 
\begin{equation*}
\widebar{\Delta} : Sp(4)^{(11)} \overset{\beta_0}\to (S^{7} \vee S^{10}) \cup e^{11} \overset{\gamma_0}\to {\wedge}^2(S^{3} \cup_{\omega} e^{7}) \hookrightarrow {\wedge}^{2}Sp(4)^{(11)},
\end{equation*}
for some $\beta_0$ and $\gamma_0$.
Then for dimensional reasons, 
$(\beta_0{\wedge}\beta_0){\comp}\alpha_0 : Sp(4)^{(15)} \to ((S^{7} \vee S^{10}) \cup e^{11}){\wedge}((S^{7} \vee S^{10}) \cup e^{11})$ and $(\gamma_0{\wedge}\gamma_0)\vert_{S^{7}{\wedge}S^{7}} : S^{7}{\wedge}S^{7} \to {\wedge}^{4}(S^{3} \cup_{\omega} e^{7})$ are respectively equal to the compositions
\begin{align*}&
(\beta_0{\wedge}\beta_0){\comp}\alpha_0 : Sp(4)^{(15)} \overset{\alpha'_0}\to S^{7}{\wedge}S^{7} \hookrightarrow ((S^{7} \vee S^{10}) \cup e^{11}){\wedge}((S^{7} \vee S^{10}) \cup e^{11}),
\\&
(\gamma_0{\wedge}\gamma_0)\vert_{S^{7}{\wedge}S^{7}} : S^{7}{\wedge}S^{7} \overset{\gamma'_0}\to {\wedge}^4S^{3} \hookrightarrow {\wedge}^{4}(S^{3} \cup_{\omega} e^{7}),
\end{align*}
for some $\alpha'_0$ and $\gamma'_0$.
Hence $\widebar{\Delta}_{4} : Sp(4)^{(15)} \to {\wedge}^{4}Sp(4)^{(15)}$ equals the composition
\begin{equation*}
\widebar{\Delta}_{4} : Sp(4)^{(15)} \overset{\alpha'_0}\to S^{7}{\wedge}S^{7} \overset{\gamma'_0}\to {\wedge}^4S^{3} \hookrightarrow {\wedge}^{4}Sp(4)^{(15)},
\end{equation*}
where $Sp(4)^{(15)} = Sp(3)^{(14)} \cup e^{15}$.
By Theorem \ref{prop:main}, $x_{7}^2 = 0$ in $h^*(Sp(3))$, and hence $\alpha'_0$ annihilates $Sp(3)^{(14)}$. 
Thus $\widebar{\Delta}_{4} : Sp(4)^{(15)} \to {\wedge}^{4}Sp(4)^{(15)}$ equals the composition 
\begin{equation*}
\widebar{\Delta}_{4} : Sp(4)^{(15)} \overset{q''}\to S^{15} \overset{\beta'_0}\to S^{14} \overset{\gamma'_0}\to S^{12} \overset{i''}\hookrightarrow {\wedge}^{4}Sp(4)^{(15)}
\end{equation*}
for some $\beta'_0$, where $q'' : Sp(4)^{(15)} \to Sp(4)^{(15)}/Sp(4)^{(14)}=S^{15}$ is the projection and $i'' : S^{12} = S^3{\wedge}S^3{\wedge}S^3{\wedge}S^3 \hookrightarrow {\wedge}^{4}Sp(4)^{(15)}$ is the inclusion.
Hence the non-triviality of $\widebar{\Delta}_{4}$ implies the non-triviality of $\beta'_0$ and $\gamma'_0$.
Therefore $\widebar{\Delta}_{4}$ should be $i''{\comp}\eta_{12}^3{\comp}q''$, if it were non-trivial.
However, we also know from (5.5) of \cite{Toda:comp-methods} that $\eta_{12}^3$ is $12\nu_{12}=6\omega_{12}$ and that $i''{\comp}\omega_{12}$ is trivial by Fact \ref{prop:nu'}.
Therefore, $\widebar{\Delta}_{4} : Sp(4)^{(15)} \to {\wedge}^{4}Sp(4)^{(15)}$ is trivial, and hence $\wcat{Sp(4)^{(15)}} \leq 3$.
This implies that $\wcat{Sp(4)^{(15)}} = 3$.
\end{Proof}
Secondly, we show the following.
\begin{Prop}\label{prop:Sp(n)-18}~
$\wcat(Sp(n)^{(19)}) = 4$ for $n \geq 4$.
\end{Prop}
\begin{Proof}~
Since $\widebar{\Delta}_{5} = ((1_{Sp(n)}){\wedge}\widebar{\Delta}_{4}){\comp}\widebar{\Delta} : Sp(n)^{(19)} \to {\wedge}^{5}Sp(n)^{(19)}$, it equals the composition
\begin{align*}&
\widebar{\Delta}_{5} : Sp(n)^{(19)} \overset{\widebar{\Delta}}\to Sp(n)^{(16)}{\wedge}Sp(n)^{(16)} = Sp(4)^{(15)}{\wedge}Sp(4)^{(15)}
\\&\qquad\qquad\qquad
\overset{(1_{Sp(4)^{(15)}}){\wedge}\widebar{\Delta}_{4}}\to {\wedge}^{5}Sp(4)^{(15)} \hookrightarrow {\wedge}^{5}Sp(n)^{(19)},
\end{align*}
which is trivial, since $\widebar{\Delta}_{4} : Sp(4)^{(15)} \to {\wedge}^{4}Sp(4)^{(15)}$ is trivial by Proposition \ref{prop:Sp(n)-15}.
Thus $\wcat(Sp(n)^{(19)}) \leq 4$, and hence $\wcat(Sp(n)^{(19)}) = 4$.
\end{Proof}

Let $p_j : Sp(n) \to X_{n,j} = Sp(n)/Sp(n-j)$ be the projection for $j \geq 1$.
Then we have the following.
\begin{Prop}\label{prop:injective-n}~
Let $q''' : Sp(n) \to Sp(n)/Sp(n)^{((2n+1)n-3)} = S^{(2n+1)n}$ be the collapsing map and $i''' : S^{(2n+1)n-6} \hookrightarrow ({\wedge}^{5}Sp(n)){\wedge}X_{n,n-3}{\wedge}\cdots{\wedge}X_{n,1}$ the inclusion.
Then 
\begin{equation*}
{q'''}^{\ast}{\comp}{i'''}_{\ast} : \pi_{(2n+1)n}(S^{(2n+1)n-6}) \to [Sp(n),({\wedge}^{5}Sp(n)){\wedge}X_{n,n-3}{\wedge}\cdots{\wedge}X_{n,1}]
\end{equation*}
is injective.
\end{Prop}
\begin{Proof}~
Firstly, we have the following exact sequence
\begin{align*}&
\pi_{(2n+1)n}(S^{(2n+1)n-3}) \overset{\psi'''}\to \pi_{(2n+1)n}(S^{(2n+1)n-6}) 
\\&\qquad
\overset{{i'''}_{\ast}}\to \pi_{(2n+1)n}(({\wedge}^{5}Sp(n)){\wedge}X_{n,n-3}{\wedge}\cdots{\wedge}X_{n,1}) \to \pi_{(2n+1)n}({\vee_5}S^{(2n+1)n-2}),
\end{align*}
where $\pi_{(2n+1)n}(S^{(2n+1)n-6}) \cong \integral/2\integral\nu^2_{(2n+1)n-6}$ and $\psi'''$ is induced from $\omega_{(2n+1)n-6}$ $=$ $2\nu_{(2n+1)n-6}$.
Thus $\psi'''$ is trivial, and hence ${i'''}_{\ast}$ is injective.

Secondly, since $({\wedge}^{5}Sp(n)){\wedge}X_{n,n-3}{\wedge}\cdots{\wedge}X_{n,1})$ is $(n(2n+1)-11)$-connected, we have 
\begin{align*}&
[Sp(n),({\wedge}^5Sp(n)){\wedge}X_{n,n-3}{\wedge}\cdots{\wedge}X_{n,1}] 
\\&\ 
= [(S^{(2n+1)n-7} \cup_{\omega_{(2n+1)n-7}} e^{(2n+1)n-3}) \vee S^{(2n+1)n},({\wedge}^5Sp(n)){\wedge}X_{n,n-3}{\wedge}\cdots{\wedge}X_{n,1}] 
\\&\ 
= [S^{(2n+1)n-7} \cup_{\omega_{(2n+1)n-7}} e^{(2n+1)n-3},({\wedge}^5Sp(n)){\wedge}X_{n,n-3}{\wedge}\cdots{\wedge}X_{n,1}]
\\&\qquad\qquad\qquad
{\oplus}\pi_{(2n+1)n}(({\wedge}^5Sp(n)){\wedge}X_{n,n-3}{\wedge}\cdots{\wedge}X_{n,1})
\end{align*}
by Proposition \ref{prop:phi}, and hence ${q'''}^{\ast}$ is injective.
Thus ${q'''}^{\ast}{\comp}{i'''}_{\ast}$ is injective.
\end{Proof}
Then the following lemma implies that $((1_{\wedge^5Sp(n)}){\wedge}p_{n-3}{\wedge}\cdots{\wedge}p_{1}){\comp}\widebar{\Delta}_{n+2}$ is non-trivial by Proposition \ref{prop:injective-n}, and hence we obtain Theorem \ref{thm:sp(n)}.
\begin{Lem}\label{lem:key-n}~
$((1_{\wedge^5Sp(n)}){\wedge}p_{n-3}{\wedge}p_{n-4}{\wedge}\cdots{\wedge}p_{1}){\comp}\widebar{\Delta}_{n+2} = {i'''}{\comp}\nu_{(2n+1)n-6}^2{\comp}{q'''}$.
\end{Lem}
\begin{Proof}~
We have 
\begin{align*}&
((1_{\wedge^5Sp(n)}){\wedge}p_{n-3}{\wedge}\cdots{\wedge}p_{1}){\comp}\widebar{\Delta}_{n+2} 
= (\widebar{\Delta}_{5}{\wedge}p_{n-3}{\wedge}\cdots{\wedge}p_{1}){\comp}\widebar{\Delta}_{n-2}
\\&\qquad
= (\widebar{\Delta}_{5}{\wedge}(1_{\wedge^{n-3}Sp(n)})){\comp}((1_{Sp(n)}){\wedge}p_{n-3}{\wedge}\cdots{\wedge}p_{1}){\comp}\widebar{\Delta}_{n-2}.
\end{align*}
For dimensional reasons, the image of $((1_{Sp(n)}){\wedge}p_{n-3}{\wedge}\cdots{\wedge}p_{1}){\comp}\widebar{\Delta}_{n-2}$ lies in 
\begin{equation*}
Sp(n)^{(21)}{\wedge}S^{15}{\wedge}\cdots{\wedge}S^{4n-1}
\cup 
Sp(n)^{(19)}{\wedge}X_{n,n-3}{\wedge}\cdots{\wedge}X_{n,1}.
\end{equation*}
From Proposition \ref{prop:Sp(n)-18}, it follows that $\widebar{\Delta}_{5}$ annihilates $Sp(n)^{(19)}$, and hence it equals the composition 
\begin{equation*}
\widebar{\Delta}_{5} : Sp(n)^{(21)} \to S^{21} \overset{\delta}\to {\wedge}^5Sp(n)^{(21)}
\end{equation*}
for some $\delta \in \pi_{21}({\wedge}^5Sp(3))$.
Then we obtain the following diagram except for the dotted arrow using Lemma \ref{lem:key}, which is commutative up to homotopy:
\begin{equation*}
\begin{xy}
\xymatrix{
Sp(n)^{(21)}
   \ar[rrrr]_{\widebar{\Delta}_{5}}
   \ar[drr]
&{}
&{}
&{}
&{\wedge}^5Sp(n)^{(21)}
\\
{}
&{}
&S^{21}
   \ar[urr]^{\delta}
   \ar@{..>}[drr]^{\delta_0}
&{}
&
\\
Sp(3)
   \ar@{^{(}->}[uu]
   \ar[rrrr]_{\widebar{\Delta}_{5}}
   \ar[dr]_{q'}
   \ar[urr]^{q'}
&{}
&{}
&{}
&{\wedge}^5Sp(3)
   \ar@{^{(}->}[uu]_{j}
\\
{}
&S^{21}
   \ar[rr]_{\nu_{15}^2}
&{}
&S^{15}
   \ar@{^{(}->}[ur]_{i'}
&{}
}
\end{xy}
\end{equation*}
Since the pair (${\wedge}^5Sp(n)$,${\wedge}^5Sp(3)$) is $26$-connected for $n \geq 4$, we can compress $\delta$ into ${\wedge}^5Sp(3)$ as $\delta \sim j{\comp}\delta_0$.
Thus we have 
$$
j{\comp}\delta_{0}{\comp}q' \sim j{\comp}i'{\comp}\nu^2_{15}{\comp}q'.
$$
Now we know that $\dim{Sp(3)} = 21 < 26-1$, and hence we can drop $j$ from the above homotopy relation and obtain 
$$
\delta_{0}{\comp}q' \sim i'{\comp}\nu^2_{15}{\comp}q'.
$$
By Proposition \ref{prop:injective}, ${q'}^{\ast} : \pi_{21}({\wedge}^5Sp(3)) \to [Sp(3),{\wedge}^5Sp(3)]$ is injective, and hence we obtain 
$$
\delta_{0} \sim i'{\comp}\nu^2_{15}.
$$
Thus $\widebar{\Delta}_{5}$ equals the composition 
\begin{equation*}
\widebar{\Delta}_{5} : Sp(n)^{(21)} \to S^{21} \overset{\nu_{15}^2}\to S^{15} \hookrightarrow {\wedge}^5Sp(n)^{(21)}.
\end{equation*}
Thus $((1_{\wedge^5Sp(n)}){\wedge}p_{n-3}{\wedge}\cdots{\wedge}p_{1}){\comp}\widebar{\Delta}_{n+2}$ equals the composition
\begin{align*}&
((1_{\wedge^5Sp(n)}){\wedge}p_{n-3}{\wedge}\cdots{\wedge}p_{1}){\comp}\widebar{\Delta}_{n+2} : Sp(n) \to S^{21}{\wedge}S^{(2n+7)(n-3)}
\\&\qquad\qquad\qquad
\overset{\nu_{(2n+1)n-6}^2}\to S^{15}{\wedge}S^{(2n+7)(n-3)} \hookrightarrow (\wedge^5Sp(n)){\wedge}X_{n,n-3}{\wedge}\cdots{\wedge}X_{n,1}.
\end{align*}
This completes the proof of the lemma.
\end{Proof}
%
%
\ifx\undefined\bysame
\newcommand{\bysame}{\leavevmode\hbox to3em{\hrulefill}\,}
\fi

\end{document}